\documentclass[11pt]{article}
\usepackage{latexsym,amsfonts}
\usepackage{amssymb,amsmath}
\usepackage{enumitem}
\usepackage{graphics}
\usepackage[demo]{graphicx}
\usepackage{epsfig}
\usepackage{pstricks}
\usepackage[normal]{subfigure}
\usepackage[latin1]{inputenc}
\usepackage[english,francais]{babel}
\usepackage{relsize,exscale}
\usepackage{makeidx}
\usepackage{enumitem}
\usepackage{amsfonts,amssymb,amsmath}
\usepackage{graphicx}
\usepackage{color}
\usepackage{multirow}
\usepackage{mathrsfs}
\usepackage[normalem]{ulem}
\usepackage{cancel}
\usepackage{bbm}
\usepackage{empheq}
\newenvironment{prooff}{{\it Proof :}}{\hfill\rule{2mm}{2mm}\vskip3mm\par}
\newtheorem{theorem}{Theorem}[section]
\newtheorem{lemma}[theorem]{Lemma}

\newtheorem{e-definition}[theorem]{Definition\rm}
\newtheorem{remark}{\it Remark\/}

\usepackage{color}
\definecolor{dred}{rgb}{0.92,0,0}
\definecolor{dgreen}{rgb}{0,0.92,0}
\definecolor{dblue}{rgb}{0,0,0.92}
\definecolor{dyellow}{rgb}{0.95,0.95,0}

\newcommand{\R}{\mathbb{R}}
\newcommand{\N}{\mathbb{N}}
\def\D{\displaystyle}
\newcommand{\ds}{\displaystyle}
\newcommand{\hs}{\hspace{0.1cm}}

\newcommand{\sa}{\\ [0.2cm]}
\usepackage{multirow}
\DeclareMathOperator*{\esssup}{ess\,sup}
\graphicspath{
{./Figures/}
{Figures/}
{./}
}
\title{Improved $P_1$-interpolation error estimates in $W^{1,p}(]0,1[)$:
Application to finite element method}
\author{Jo\"el Chaskalovic \thanks{D'Alembert,
Sorbonne University, Paris, France, (Email:  jch1826@gmail.com)}
\qquad
Franck Assous
\thanks{
Department of Mathematics, Ariel University, Ariel, Israel, (Email: assous@ariel.ac.il)}
\qquad
}
\date{}
\begin{document}
\maketitle
\selectlanguage{english}
\begin{abstract}
\noindent Based on a new Taylor-like formula, we derived an improved interpolation error estimate in $W^{1,p}$. We compare it with the classical error estimates based on the standard Taylor formula, and also with the corresponding interpolation error estimate, derived from the mean value theorem. We then assess the improvement in accuracy we can get from this formula, leading to a significant reduction in finite element computation costs.
\noindent
\end{abstract}
\noindent {\em keywords}: Taylor's theorem, Taylor-like formula, Error estimate, Interpolation error, Approximation error, Finite elements.
\section{Introduction}\label{intro}

\noindent Even today, improving the precision of approximations continues to be a challenge in the field of numerical analysis. In this context, we have recently introduced a
second-order Taylor-like formula \cite{JCAM2023}, which is itself an extension of a first-order Taylor-like formula published earlier (see \cite{arXiv_First_Order}).
The principle is as follows: if we view the classical Taylor formula as consisting of two parts - the polynomial component and the remainder - the main idea behind these
formulas is to minimize the remainder by redistributing some of its numerical weight  to the polynomial component.\sa
To achieve this, we introduced a sequence of $(n+1)$ equally spaced points and examined a linear combination of the first derivatives at these points. We then demonstrated
that an optimal selection of the coefficients in this linear combination minimizes the associated remainder, particularly when compared to that of the standard
Taylor's formula \cite{Taylor}. One of the major consequences of these new formulas is obtaining smaller constants in error estimates which play a significant role in
assessing the accuracy of numerical approximation methods.\sa
Indeed, in all error estimates based on Taylor's formula, there are constants that cannot be precisely calculated but can only be bounded. This is due to the existence of a
non-unique unknown point that appears in the remainder of Taylor's expansion, inherited from Rolle's theorem. Consequently, the accuracy of a given numerical method is often
assessed by examining the asymptotic convergence rate as the mesh size involved in the numerical method tends to zero. However, this situation does not correspond to any
concrete application for which the mesh size is fixed.\sa
So, let us consider, for instance, approximation errors in the finite element method where $u$ represents the exact solution to a variational problem, and $u_h^{(k)}$ and $u_h^{(m)}, (k<m),$ denote the corresponding $P_k$ and $P_m$ finite element approximations. \sa
In this context, the error estimates have the following structure (see for instance \cite{ChaskaPDE},
\cite{Ern_Guermond}, \cite{BrSc08} and \cite{RaTho82}):
$$
\D \|u-u_h^{(k)}\| \leq C_k\,h^k, \, \mbox{ and } \, \|u-u_h^{(m)}\| \leq C_m\,h^m,
$$
for a given norm $\|.\|$ which can be specified, depending on the concerned application. \sa
Here, the constants $C_k$ and $C_m$ are unknown, amongst other, stemming from the uncertainty related to the unknown point of Taylor's formula \cite{AsCh2014}.\sa
However, even for a small, yet {\em fixed value} of the mesh size $h$, because the unknown value of $C_k$ and $C_m$, it can occur that $C_k\,h^k \leq C_m\,h^m$. Consequently,
one cannot determine whether the $P_m$ finite element method is more accurate than the $P_k$ one, as we lack information about the relative positions between $\|u-u_h^{(k)}\|$
and $\|u-u_h^{(m)}\|$, which is the main issue, even if $C_m\,h^m \leq C_k\,h^k$. \sa
For this reason, asymptotic behavior is generally considered, because as $h$ tends to 0, $h^m$ converges to 0 more rapidly than $h^k$. This results in $\|u-u_h^{(m)}\|$
converging to zero faster than $\|u-u_h^{(k)}\|$, indicating that the $P_m$ finite elements are more accurate than the $P_k$ ones {\em in an asymptotic sense}. However,
in practical applications, the mesh size is fixed and does not tend to zero, and the asymptotic comparison is then not valid anymore. \sa
That is why we focus our attention towards improving the evaluation of the constants $C_k$ and $C_m$.
We ask ourselves if the constants, typically usually in error estimates, are as small as possible.\sa
In this context, various approaches have been proposed to enhance the approximation accuracy. For instance, in the field of numerical integration, readers can refer
to \cite{Barnett_Dragomir}, \cite{Cerone} or \cite{Dragomir_Sofo}, along with the references cited therein. From a different perspective, due to the lack of information,
heuristic methods have been explored, primarily based on a probabilistic approach, see for instance \cite{Abdulle}, \cite{AsCh2014}, \cite{Hennig}, \cite{Oates}
or \cite{ArXiv_JCH}, \cite{CMAM2} and \cite{ChAs20}. These approaches enable the comparison of different numerical methods for a fixed mesh size as shown in \cite{MMA2021}.\sa
We previously published several improved error estimates in \cite{arXiv_First_Order} and \cite{JCAM2023}. In this paper, we continue this exploration by examining the influence
of the Taylor-like formula on the $W^{1,p}$ interpolation error estimate in the one-dimensional case. To achieve this, we consider three different methods: the standard
Taylor formula,  the mean value theorem, and  the new Taylor-like formula.\sa
We compare the constants in these three error estimates and assess the benefits within the context of finite element applications. In particular, we show a significant reduction
in the number of nodes required, thanks to the reduction of the constant obtained, using the Taylor-like formula.\sa
The paper is organized as follows. In Section \ref{Interpolation}, we present the primary results, focusing on $W^{1,p}$ interpolation error estimates. Section \ref{Finite_Element} explores the implications of these interpolation error estimates within the framework of finite elements. We also provide examples in various $W^{m,p}$ Sobolev spaces to illustrate the new interpolation error estimates and their impact on implementation efficiency. Finally, we conclude with some remarks.
\section{Improving $P_1$-interpolation error estimate in $W^{1,p}(]0,1[)$}\label{Interpolation}
\noindent In this section, we derive a new $W^{1,p}-$ interpolation error estimate, based on the Taylor-like formula derived in \cite{arXiv_First_Order}. Since the
special case $p=1$ has been addressed in \cite{ChAs23W11}, we assume in the sequel of the paper that $p$ is an integer strictly greater than one.\sa
We consider a given real function $u$ defined on the interval $[0,1]$ which belongs to $C^2([0,1])$. Then, $\exists\, (m_2,M_2)\in\R^2$ such that,
$\forall x \in [0,1]: m_{2} \leqslant u''(x) \leqslant M_{2}$.\sa
We also introduce a sequence of $N+2$ points $(x_i)_{i=0,\dots,N+1}$ in $[0,1]$ defined by
\begin{empheq}[left=\hspace{-11cm}\empheqlbrace]{alignat=2}
\D \hs  & x_0  & \hs = \hs\hs & 0,\, x_{N+1}=1, \nonumber \\[0.2cm]
\D \hs & x_{i+1} & \hs = \hs\hs & x_i + h_i, \, (i=0,\dots,N), \nonumber
\end{empheq}
and we define the mesh size $h$ by: $\D h=\max_{i=0,\dots,N}h_i$. \sa
Finally, we consider the $P_1$-interpolation polynomial $u_{I}$ of $u$ which satisfies
\begin{eqnarray}
\forall i \in \{0,\dots,N+1\}, \, u_{I}(x_i)=u(x_i), & & \nonumber\\[0.2cm]
\forall x \in [x_i,x_{i+1}], u_{I} \in P_1([x_i,x_{i+1}]), & & \nonumber
\end{eqnarray}
where $P_1([x_i,x_{i+1}])$ denotes the space of polynomials defined on $[x_i,x_{i+1}]$ of degree less than or equal to 1.\sa
In the following, we also need some notations. We consider the Sobolev space $W^{1,p}(]0,1[)$ defined by
$$
\D W^{1,p}(]0,1[) = \left\{u:\,]0,1[\rightarrow\R, \frac{}{}u \in L^p(]0,1[)\,;\, u' \in L^p(]0,1[)\right\},
$$
where $u'$ denotes the weak derivative of $u$ which belongs to $L^p(]0,1[)$, \cite{Brezis}. For any $u\in W^{1,p}(]0,1[)$,
we denote by $\|.\|_{1,p}$ the classical norm defined by
$$
\D \|u\|_{1,p} = \bigg(\|u\|^p_{0,p} + \|u'\|^p_{0,p}\bigg)^{\!\frac{1}{p}},
$$
the norm $\|.\|_{0,p}$ being defined by
$$
\D \forall u \in L^p(]0,1[): \|u\|_{0,p} = \bigg[\int_{0}^{1}|u(x)|^p\,dx\bigg]^{\frac{1}{p}}.
$$
We will first derive a useful lemma that will be applied several times later in this paper.
\begin{lemma}\label{Lem-Generic}
Let $u\in C^2([0,1])$. $\forall i=0,\dots,N$, we set:
\begin{equation}\label{x'k}
\D x'_k=x_i+\frac{kh_i}{n}, \forall k=0,\dots,n,
\end{equation}
where $n$ is a given non-zero integer. \sa
Then
$$
\D \int_{x_i}^{x_{i+1}}|u'(x)-u'(x'_k)|^p\, dx \leq \frac{1}{p+1}\bigg[k^{p+1} + (n-k)^{p+1}\bigg]\bigg(\frac{h_i}{n}\bigg)^{p+1}\|u''\|^p_{\infty},
$$
where $\D\|u''\|_{\infty}=\esssup_{x\in[0, 1]}|u''(x)|.$
\end{lemma}
\begin{prooff}
Let $x'_k\in[x_i,x_{i+1}]$ defined by (\ref{x'k}). We have:
\begin{eqnarray}
\D \int_{x_i}^{x_{i+1}}|u'(x)-u'(x'_k)|^p\, dx & = & \int_{x_i}^{x_{i+1}}\bigg|\int_{x'_k}^{x}u''(t)\,dt\bigg|^p\, dx, \nonumber\\[0.2cm]
\D & = & \int_{x_i}^{x'_{k}}\bigg|\int_{x'_k}^{x}u''(t)\,dt\bigg|^p\, dx + \int_{x'_{k}}^{x_{i+1}}\bigg|\int_{x'_k}^{x}u''(t)\,dt\bigg|^p\, dx, \nonumber\\[0.2cm]
\D & \leq & \int_{x_i}^{x'_{k}}\bigg(\int_{x}^{x'_k}|u''(t)|\,dt\bigg)^p\, dx + \int_{x'_{k}}^{x_{i+1}}\bigg(\int_{x'_k}^{x}|u''(t)|\,dt\bigg)^p dx. \label{II1}
\end{eqnarray}
Let us now introduce the non-zero integer $q$, the conjugate of $p$, which satisfies
\begin{equation}\label{pq}
\D \frac{1}{p}+\frac{1}{q}=1.
\end{equation}
Then, by the help of Hölder's inequality \cite{Brezis}, (\ref{II1}) leads to
\begin{equation}\label{II2}
\D \int_{x_i}^{x_{i+1}}|u'(x)-u'(x'_k)|^p\, dx \leq \int_{x_i}^{x'_{k}}(x'_k-x)^{\frac{p}{q}}\bigg(\int_{x}^{x'_k}|u''(t)|^p dt\bigg) dx \, + \int_{x'_{k}}^{x_{i+1}}(x-x'_k)^{\frac{p}{q}}\bigg(\int_{x'_k}^{x}|u''(t)|^p dt\bigg) dx.
\end{equation}
Let us remark that given (\ref{pq}), $\D\frac{p}{q}=p-1$. So, the function $F$ defined by:
$$F(x,t)=(x'_k-x)^{p-1}|u''(t)|^p, \forall (x,t)\in [0,1]\times [0,1],$$
fulfills Tonelli's theorem \cite{Brezis}, for all $p\geq 1$. As a consequence, $F$ belongs to $L^1\!\big([0,1]\times [0,1]\big)$, and is {\em a fortiori} in
$L_{loc}^1\!\big([0,1]\times [0,1]\big)$. \sa
Therefore, Fubini's theorem \cite{Brezis} can be applied to the first integral of (\ref{II2}). Since the same arguments are valid for the second integral of (\ref{II2}), we have
$$
\D \int_{x_i}^{x_{i+1}}\!\!|u'(x)-u'(x'_k)|^p\, dx \leq \int_{x_i}^{x'_{k}}|u''(t)|^p\bigg(\!\!\int_{x_i}^{t}(x'_k-x)^{p-1}\,dx\!\!\bigg) dt\, + \int_{x'_{k}}^{x_{i+1}}\!\!|u''(t)|^p\bigg(\!\!\int_{t}^{x_{i+1}}\!\!(x-x'_k)^{p-1}dx\!\!\bigg) dt,
$$
which can be written, after integrating which respect to $t$:
\begin{eqnarray}
\D \int_{x_i}^{x_{i+1}}|u'(x)-u'(x'_k)|^p\, dx \, \leq\, & \hspace{-0.2cm}\D\frac{1}{p}&\hspace{-0.3cm}\bigg(
\int_{x_i}^{x'_{k}}|u''(t)|^p\big[(x'_k-x_i)^{p}-(x'_k-t)^{p}\big] dt \nonumber\\[0.2cm]
\D & \hspace{-0.2cm}+ & \hspace{-0.3cm}\int_{x'_{k}}^{x_{i+1}}|u''(t)|^p \big[(x_{i+1}-x'_k)^{p}-(t-x'_k)^{p}\big]dt\bigg).\label{II4}
\end{eqnarray}
Regarding the first integral of the right-hand side of (\ref{II4}), $t$ belonging to the interval $[x_i,x'_k]$, we set:
$$t=sx_i + (1-s)x'_k, \, (s\in[0,1]),$$
and we get that
\begin{eqnarray}
\D \int_{x_i}^{x'_{k}}|u''(t)|^p\big[(x'_k-x_i)^{p}-(x'_k-t)^{p}\big] dt & \leq & \|u''\|^p_{\infty}\int_{0}^{1}(x'_k-x_i)^{p+1}(1-s^p) ds, \nonumber\\[0.2cm]
& \leq & \frac{p}{p+1}(x'_k-x_i)^{p+1}\|u''\|^p_{\infty}. \label{II5}
\end{eqnarray}
In the same way, we obtain for the second integral of the right-hand side of (\ref{II4}):
\begin{equation}\label{II6}
\D \int_{x'_{k}}^{x_{i+1}}|u''(t)|^p \big[(x_{i+1}-x'_k)^{p}-(t-x'_k)^{p}\big] dt \leq \frac{p}{p+1}(x_{i+1}-x'_k)^{p+1}\|u''\|^p_{\infty}.
\end{equation}
Summing up (\ref{II5}) and (\ref{II6}) and dividing by $p$, inequality (\ref{II4}) gives
\begin{equation}\label{IIhelp}
\D \int_{x_i}^{x_{i+1}}|u'(x)-u'(x'_k)|^p\, dx \leq  \frac{1}{p+1}\bigg[(x'_k-x_i)^{p+1} + (x_{i+1}-x'_k)^{p+1}\bigg]\|u''\|^p_{\infty}.
\end{equation}
Now considering that $\D x'_k-x_i=\frac{kh_i}{n}$ and $\D x_{i+1}-x'_k = \bigg(1-\frac{k}{n}\bigg)h_i$, this leads to
$$
\D \int_{x_i}^{x_{i+1}}|u'(x)-u'(x'_k)|^p\, dx \leq \frac{1}{p+1}\bigg[k^{p+1} + (n-k)^{p+1}\bigg]\bigg(\frac{h_i}{n}\bigg)^{p+1}\|u''\|^p_{\infty}.
$$
\end{prooff}
In the first step, we will derive the interpolation error estimate based on the standard Taylor formula. Subsequently, we will derive the analogous result obtained using the Taylor-like formula. Let us begin with the following lemma:
\begin{lemma}\label{Error_Taylor_Lem}
Let u be in $C^2([0,1])$ and $u_I$ the corresponding $P_1$-interpolation polynomial. Then, the standard Taylor formula leads to the interpolation error estimate:
\begin{equation}\label{Error_Estim_T}
\D \|u-u_I\|^p_{1,p} \leq \bigg(\frac{2^{p-1}}{p+1} +\frac{1}{2}\bigg)h^p\bigg(1+\frac{h^p}{p}\bigg)\|u''\|^p_{\infty}.
\end{equation}
\end{lemma}
\begin{prooff}
We recall the classical first order Taylor formula as expressed in \cite{arXiv_First_Order}:
\begin{equation}\label{Taylor}
u(x_{i+1}) = u(x_i) + h_i u'(x_i) + h_i \epsilon^{(T)},
\end{equation}
with
\begin{equation}\label{Reste_Taylor}
\D|\epsilon^{(T)}| \leq \frac{h_i}{2}\|u''\|_{\infty}.
\end{equation}
We begin by evaluating the $L^{p}$-norm of the derivative, that is $\|u'-u'_{I}\|_{0,p}$.\sa
We have:
$$
\D \|u'-u'_{I}\|^p_{0,p} = \int_{0}^{1}|u'(x)-u'_{I}(x)|^p\,dx = \sum_{i=0}^{N}\int_{x_i}^{x_{i+1}}|u'(x)-u'_{I}(x)|^p\,dx.
$$
Then, given that $u'_I$ is constant on $[x_i,x_{i+1}]$,  by the help of (\ref{Taylor}), we have
$$
\D \forall x \in [x_i,x_{i+1}]: u'_{I}(x) = \frac{u(x_{i+1})-u(x_i)}{h_i} = u'(x_i) + \epsilon^{(T)}.
$$
As a consequence
\begin{eqnarray}
\D \int_{x_i}^{x_{i+1}}|u'(x)-u'_{I}(x)|^p\,dx & = & \int_{x_i}^{x_{i+1}}|u'(x)-u'(x_i) - \epsilon^{(T)}|^p\,dx,
\nonumber\\[0.2cm]
\D & \leq & 2^{\frac{p}{q}}\int_{x_i}^{x_{i+1}}|u'(x)-u'(x_i)|^p\, dx +  \frac{h_i^{p+1}}{2}\|u''\|_{\infty}^{p}, \label{IN2}
\end{eqnarray}
where $p$ and $q$ are conjugates as defined in (\ref{pq}), the reminder $\epsilon^{(T)}$ is bounded from (\ref{Reste_Taylor}). Furthermore, to obtain (\ref{IN2}), we
also used the discrete Hölder's inequality \cite{Brezis} in the following way:
$$
\D \forall \, i=1,\dots,m,\,  \forall (a_i,b_i)\in\R^2: \D\sum_{i=1}^{m}|a_ib_i| \leq \bigg(\D\sum_{i=1}^{m}|a_i|^p\bigg)^{\frac{1}{p}}\bigg(\D\sum_{i=1}^{m}|b_i|^q\bigg)^{\frac{1}{q}}.
$$
Then, if $b_i=1$, for all $i=1,\dots,m$, this inequality becomes
\begin{equation}\label{Holder_part1}
\D\sum_{i=1}^{m}|a_i| \leq m^{\frac{1}{q}}\bigg(\D\sum_{i=1}^{m}|a_i|^p\bigg)^{\frac{1}{p}}.
\end{equation}
Specifically, when $m=2$, we obtain that
$$
\D |a_1+a_2|^p \leq 2^{\frac{p}{q}}\bigg(|a_1|^p+|a_2|^p\bigg)\,,
$$
that gives (\ref{IN2}) by choosing $a_1=u'(x)-u'(x_i)$ and $a_2=\epsilon^{(T)}$.\sa
\sa
Let us apply Lemma \ref{Lem-Generic} to the integral on the right-hand side of (\ref{IN2}) by choosing the point $x'_k=x_i$. This leads to
$$
\D \int_{x_i}^{x_{i+1}}|u'(x)-u'(x_i)|^p\,dx \leq \frac{1}{p+1}h_i^{p+1}\|u''\|^p_{\infty},
$$
and (\ref{IN2}) becomes
\begin{equation}\label{IN4}
\D \int_{x_i}^{x_{i+1}}|u'(x)-u'_{I}(x)|^p\,dx \leq \bigg(\frac{2^{p-1}}{p+1} +\frac{1}{2}\bigg) h_i^{p+1}\|u''\|^p_{\infty}.
\end{equation}
Finally, by summing in (\ref{IN4}) over $i$ from 0 to $N$, and using that $\D\sum_{i=0}^{N}h_i=1$, as well as  $h_i\leq h$, we get that
\begin{equation}\label{N1_u1}
\D \|u'-u'_{I}\|^p_{0,p} \leq \bigg(\frac{2^{p-1}}{p+1} +\frac{1}{2}\bigg)\bigg(\sum_{i=0}^{N}h_i^{p+1}\bigg)\|u''\|^p_{\infty} = \bigg(\frac{2^{p-1}}{p+1} +\frac{1}{2}\bigg)h^{p}\|u''\|^p_{\infty}\,.
\end{equation}
\sa
Let us now evaluate the $L^{p}$-norm $\|u-u_{I}\|_{0,p}$. To begin, recalling that $u(x_i)=u_I(x_i)$, we have, for all $x \in [x_i,x_{i+1}]$
$$
|u(x)-u_{I}(x)|^p = \bigg| \int_{x_i}^{x}\big(u'(t)-u_I'(t)\big)\,dt \bigg|^p.
$$
Now, using Hölder's inequality, we can write
\begin{equation}\label{Ineq00}
\D |u(x)-u_{I}(x)|^p \leq (x-x_i)^{\frac{p}{q}}\int_{x_i}^{x}|u'(t)-u'_I(t)|^p\,dt \leq (x-x_i)^{p-1}\int_{x_i}^{x_{i+1}}|u'(t)-u'_I(t)|^p\,dt.
\end{equation}
So, using inequality (\ref{IN4}), (\ref{Ineq00}) gives
$$
\D |u(x)-u_{I}(x)|^p \leq (x-x_i)^{p-1}\bigg(\frac{2^{p-1}}{p+1} +\frac{1}{2}\bigg) h_i^{p+1}\|u''\|^p_{\infty}\,.
$$
It remains now to integrate this inequality on $[x_i,x_{i+1}]$ to obtain that
$$
\D \int_{x_i}^{x_{i+1}}|u(x)-u_{I}(x)|^p dx \leq \bigg(\frac{2^{p-1}}{p+1} +\frac{1}{2}\bigg) \frac{h_i^{2p+1}}{p}\|u''\|^p_{\infty}\,,
$$
and summing over all values of $i$ between $0$ and $N$, this implies that
\begin{equation}\label{N1_u2}
\D \|u-u_{I}\|^p_{0,1} \leq \bigg(\frac{2^{p-1}}{p+1} +\frac{1}{2}\bigg) \frac{h^{2p}}{p}\|u''\|^p_{\infty}.
\end{equation}
Finally, by combining inequalities (\ref{N1_u1}) and (\ref{N1_u2}), we get the interpolation error estimate (\ref{Error_Estim_T}).
\end{prooff}
The next step consists to derive the interpolation error estimate analogous to (\ref{Error_Estim_T}), which can be obtained by using the Taylor-like formula presented in
\cite{arXiv_First_Order}. Before that, let us provide an additional result obtained by substituting the classical Taylor formula with the mean value theorem, (see for
example \cite{Atki88} or \cite{BuFa11}).
\begin{lemma}
Let u be in $C^2([0,1])$ and $u_I$ the corresponding $P_1$-interpolation polynomial. Then, the mean value theorem leads to the following interpolation error estimate:
\begin{equation}\label{Error_Estim_ACF}
\D \|u-u_I\|^p_{1,p} \leq \frac{1}{p+1}\,h^p\bigg(1+\frac{h^p}{p}\bigg)\|u''\|^p_{\infty}.
\end{equation}
\end{lemma}
\begin{prooff} Here also, we begin by evaluating the $L^p$-norm of the derivative, that is $\|u'-u'_{I}\|_{0,p}$. \sa
Given that $u'_I$ is constant on $[x_i,x_{i+1}]$,  the mean value theorem enables to write that, for all $x \in [x_i,x_{i+1}]$, there exists a point $\xi_i$ belonging to $]x_i,x_{i+1}[$ such that
$$
u'_{I}(x) = \frac{u(x_{i+1})-u(x_i)}{h_i} = u'(\xi_i)\,.
$$
Hence, we obtain that
\begin{equation}\label{IN1_2}
\D \int_{x_i}^{x_{i+1}}|u'(x)-u'_{I}(x)|^p\,dx = \int_{x_i}^{x_{i+1}}|u'(x)-u'(\xi_i)|^p\,dx.
\end{equation}
Now, following the same method as described in Lemma \ref{Lem-Generic}, one can prove that, as shown in formula (\ref{IIhelp}) by simply replacing $x'_{k}$ by $\xi_i$:
$$
\D \int_{x_i}^{x_{i+1}}|u'(x)-u'(\xi_i)|^p\,dx \leq \frac{1}{p+1}\bigg[(\xi_i-x_i)^{p+1}+(x_{i+1}-\xi_i)^{p+1}\bigg]\|u''\|^p_{\infty}.
$$
In addition, using that $\xi_i \in ]x_i, x_{i+1}[$, $\xi_i$ can be written as a convex combination of $x_i$ and $x_{i+1}$, namely,
$$
\xi_i= sx_i+(1-s)x_{i+1},\, (0<s<1)\,.
$$
In these conditions, we obtain that
\begin{equation}\label{Eq-21}
\D \int_{x_i}^{x_{i+1}}|u'(x)-u'(\xi_i)|^p\,dx \leq \frac{1}{p+1}\bigg[(1-s)^{p+1}+s^{p+1}\bigg]h_i^{p+1}\|u''\|^p_{\infty} \leq \frac{1}{p+1}h_i^{p+1}\|u''\|^p_{\infty}\,,
\end{equation}
where we used that $0 \leq (1-s)^{p+1} \leq 1-s$ and $0 \leq s^{p+1} \leq s$,  for $0 \leq s \leq 1$.
\sa
Finally, by summing over $i$ from $0$ to $N$, and still using that $\D \sum_{i=0}^{N} h_i =1$, we get that
\begin{equation}\label{ACF_1}
\D \|u'-u'_I\|^p_{0,p} \leq \frac{1}{p+1}\,h^{p}\|u''\|^p_{\infty}\,.
\end{equation}
To evaluate now the norm $ \|u-u_I\|^p_{1,p}$, we follow the same procedure as detailed in Lemma \ref{Error_Taylor_Lem}. Writing also that
$$
|u(x)-u_{I}(x)|^p = \bigg| \int_{x_i}^{x}\big(u'(t)-u_I'(t)\big)\,dt \bigg|^p\,,
$$
we use Hölder's inequality together with inequality (\ref{Eq-21}) to obtain
$$
\D |u(x)-u_{I}(x)|^p \leq (x-x_i)^{p-1}\frac{1}{p+1} h_i^{p+1}\|u''\|^p_{\infty}\,.
$$
Integrating this inequality on $[x_i,x_{i+1}]$ and summing over all values of $i$ between $0$ and $N$, this leads to
$$
\D \|u-u_{I}\|^p_{0,1} \leq \frac{1}{p+1}\frac{h^{2p}}{p}\|u''\|^p_{\infty}.
$$
So, adding this expression to (\ref{ACF_1}), we finally obtain that
\begin{equation}\label{Error_Estim_ACF_V2}
\D \|u-u_I\|^p_{1,p} \leq \frac{1}{p+1}\,h^p\bigg(1+\frac{h^p}{p}\bigg)\|u''\|^p_{\infty}.
\end{equation}
\end{prooff}
\noindent Now, let us derive the interpolation error estimate obtained by using the Taylor-like formula proposed in \cite{arXiv_First_Order}, instead of the standard Taylor formula.\sa
To this end, let us first choose an integer $n\in\N^{*}$. Then, for any function $u\in C^2([0,1])$, the Taylor-like formula can be written as
\begin{equation}\label{Generalized_Taylor}
u(x_{i+1}) = u(x_i) + h_i\left(\frac{u'(x_i) + u'(x_{i+1})}{2n} + \frac{1}{n}\sum \limits_{k=1}^{n-1} u'(x'_k)\right) + h_i\epsilon_{n},
\end{equation}
where $x'_k$ is defined in (\ref{x'k}), and the remainder $\epsilon_{n}$ is bounded by
\begin{equation}\label{majoration_epsilon_n+1}
\D|\epsilon_{n}| \leqslant \frac{h_i}{8n}(M_2-m_2)\,.
\end{equation}
%
Then, we can prove the following result:
\begin{theorem}
Let u be in $C^2([0,1])$ and $u_I$ the corresponding $P_1$-interpolation polynomial. Then, the Taylor-like formula (\ref{Generalized_Taylor}) leads to the interpolation error estimate
\begin{equation}\label{Error_Estim_TL}
\D \|u-u_{I}\|^p_{1,p} \leq \frac{(n+2)^{p-1}}{p+1}\biggl(\!\frac{1}{2^{p-1}n^p} + \frac{2S^{*}_{p}(n)}{n^{2p+1}}\!\biggr)\bigg(h^p + \frac{h^{2p}}{p}\bigg)\|u''\|^p_{\infty} + \frac{1}{3n}\bigg(\frac{3}{8}\bigg)^{\!p}\bigg(h^p + \frac{h^{2p}}{p}\bigg)(M_2-m_2)^p\!,
\end{equation}
where we set
\begin{equation}\label{Spn}
\D S^{*}_{p}(n) = \sum \limits_{k=1}^{n-1}k^{p+1} \mbox{ for } n \geq 2, \mbox{ and } S^{*}_{p}(1) = 0\,.
\end{equation}
\end{theorem}
\begin{prooff}
Here also, we begin by evaluating the $L^{p}$-norm of the derivative, that is $\|u'-u'_{I}\|_{0,p}$. By the help of (\ref{Generalized_Taylor})-(\ref{majoration_epsilon_n+1}),
we obtain that
\begin{equation}\label{Ineq03}
\int_{x_i}^{x_{i+1}}|u'(x)-u'_{I}(x)|^p\,dx = \int_{x_i}^{x_{i+1}}\bigg|u'(x) -\bigg( \frac{u'(x_i) + u'(x_{i+1})}{2n} + \frac{1}{n}\sum \limits_{k=1}^{n-1} u'(x'_k) + \epsilon_{n} \bigg) \bigg|^p\, dx.
\end{equation}
Writing now $u'(x)$ in the integral as
$$
\D u'(x) = \frac{1}{2n} u'(x) + \frac{1}{2n} u'(x) + \frac{1}{n}\sum_{k=1}^{n-1}u'(x)
$$
enables us to derive the following estimate from (\ref{Ineq03}):
\begin{equation}\label{Ineq03_V1}
\int_{x_i}^{x_{i+1}}|u'(x)-u'_{I}(x)|^p\,dx \leq \int_{x_i}^{x_{i+1}}\bigg(\frac{1}{2n}|u'(x)-u'(x_i)| +\frac{1}{2n}|u'(x)-u'(x_{i+1})|+\frac{1}{n}\sum \limits_{k=1}^{n-1} |u'(x)-u'(x'_k)|
+ |\epsilon_{n}| \bigg)^p\, dx\,.
\end{equation}
Now, considering the sum of the $n+2$ terms in the parenthesis, we use the particular case of Hölder's inequality (\ref{Holder_part1}). Hence, we obtain for the terms inside the integral, using still that $\ds\frac{p}{q}=p-1$:
\begin{eqnarray*}
&&\hspace*{-0.8cm}\bigg(
\frac{1}{2n}|u'(x)-u'(x_i)| +\frac{1}{2n}|u'(x)-u'(x_{i+1})|+\frac{1}{n}\sum \limits_{k=1}^{n-1} |u'(x)-u'(x'_k)|
+ |\epsilon_{n}|
\bigg)^p \\
& \leq &(n+2)^{p-1}\bigg(
 \frac{1}{(2n)^p}|u'(x)-u'(x_i)|^{p} +\frac{1}{(2n)^p}|u'(x)-u'(x_{i+1})|^{p}+\frac{1}{n^{p}}\sum \limits_{k=1}^{n-1} |u'(x)-u'(x'_k)|^{p}+ |\epsilon_{n}|^{p}
\bigg)\,.
\end{eqnarray*}
Consequently, inequality (\ref{Ineq03_V1}) becomes
\begin{equation}\label{I1I2I3}
\D \int_{x_i}^{x_{i+1}}|u'(x)-u'_{I}(x)|^p dx \leq (n+2)^{p-1}\bigl(I_1 + I_2 + I_3\bigr) + \frac{(n+2)^{p-1}}{(8n)^p}h_i^{p+1}(M_2-m_2)^p,
\end{equation}
where we set
\begin{eqnarray*}
\D I_1 & = & \frac{1}{(2n)^p}\int_{x_i}^{x_{i+1}}|u'(x)-u'(x_i)|^p dx, \\[0.2cm]
\D I_2 & = & \frac{1}{(2n)^p}\int_{x_i}^{x_{i+1}}|u'(x)-u'(x_{i+1})|^p dx, \\[0.2cm]
\D I_3 & = & \frac{1}{n^p}\sum \limits_{k=1}^{n-1} \int_{x_i}^{x_{i+1}}|u'(x)-u'(x'_k)|^p dx.
\end{eqnarray*}
Applying now Lemma \ref{Lem-Generic}, we derive for the integrals in $I_1$ and $I_2$ the same estimate, that is
\begin{eqnarray*}
\D \hspace{-0.2cm} I_1 &
\leq \ds\frac{1}{(p+1)(2n)^{p}}\,\bigg(\frac{k^{p+1}+(n-k)^{p+1}}{n^{p+1}}\bigg)\,h_i^{p+1}\|u''\|^p_{\infty}
&  \leq \frac{1}{(p+1)(2n)^{p}}\,h_i^{p+1}\|u''\|^p_{\infty},\\[0.2cm]
\D \hspace{-0.2cm}I_2 &
\leq \ds\frac{1}{(p+1)(2n)^{p}}\,\bigg(\frac{k^{p+1}+(n-k)^{p+1}}{n^{p+1}}\bigg)\,h_i^{p+1}\|u''\|^p_{\infty}
 & \leq  \frac{1}{(p+1)(2n)^{p}}\,h_i^{p+1}\|u''\|^p_{\infty},
\end{eqnarray*}
where we used, as in (\ref{Eq-21}), that $0 \leq \tilde{k}^{p+1} \leq \tilde{k}$ and $0 \leq (1-\tilde{k})^{p+1} \leq 1-\tilde{k}$ for $0 \leq \tilde{k}:=\ds\frac{k}{n} \leq 1$, so that
$$\ds\frac{k^{p+1}+(n-k)^{p+1}}{n^{p+1}} \leq 1.$$
Considering now the integral involved in $I_3$, we obtain that
\begin{eqnarray*}
\D \hspace{-0.2cm}I_3 & \leq & \frac{1}{(p+1)n^p}\biggl(\sum \limits_{k=1}^{n-1} \bigl[k^{p+1}+(n-k)^{p+1}\bigr]\!\biggr)\!\biggl(\frac{h_i}{n}\bigg)^{p+1}\!\!\!\|u''\|^p_{\infty} = \frac{2S^{*}_{p}(n)}{(p+1)n^{2p+1}}
h_{i}^{p+1} \|u''\|^p_{\infty},
\end{eqnarray*}
where the term $S^{*}_{p}(n)$ is defined in (\ref{Spn}).\sa
Now, combining the inequalities obtained for $I_1, I_2$ and $I_3$, and using for the last term of (\ref{I1I2I3}) that
$$\ds\frac{(n+2)^{p-1}}{(8n)^{p}} \leq \frac{(3n)^{p-1}}{(8n)^{p}} = \frac{1}{3n}\bigg(\frac{3}{8}\bigg)^{\!p}, \mbox{ for } n \geq 1,$$
we find that
\begin{equation}\label{I1I2I3_01}
\D \int_{x_i}^{x_{i+1}}\!\!|u'(x)-u'_{I}(x)|^p dx \leq \frac{(n+2)^{p-1}}{p+1}\biggl(\!\frac{1}{2^{p-1}n^p} + \frac{2S^{*}_{p}(n)}{n^{2p+1}}\!\biggr)h_i^{p+1}\|u''\|^p_{\infty}
+\frac{1}{3n}\bigg(\frac{3}{8}\bigg)^{\!p}
h_i^{p+1}(M_2-m_2)^p\,.
\end{equation}
Finally, by summing over $i$ between $0$ and $N$ in (\ref{I1I2I3_01}), and using that $\ds\sum_{i=0}^{N}h_i=1$, we get
\begin{equation}\label{N1_u1_V2}
\D \|u'-u'_{I}\|^p_{0,p} \leq \frac{(n+2)^{p-1}}{p+1}\biggl(\!\frac{1}{2^{p-1}n^p} + \frac{2S^{*}_{p}(n)}{n^{2p+1}}\!\biggr)h^{p}\|u''\|^p_{\infty}
+ \frac{1}{3n}\bigg(\frac{3}{8}\bigg)^{\!p}
h^{p}(M_2-m_2)^p.
\end{equation}
\sa
\noindent Let us now evaluate the $L^p$-norm $\|u-u_{I}\|_{0,p}$. Like in (\ref{Ineq00}) and using (\ref{I1I2I3_01}), we have, $\forall x \in [x_i,x_{i+1}],$
\begin{eqnarray*}
\D \hspace{-0.4cm}|u(x)-u_{I}(x)|^p & \!\leq \!& \!\!(x\!-\!x_i)^{p-1}\!\!\int_{x_i}^{x_{i+1}}\!\!|u'(t)-u'_I(t)|^p dt, \\[0.2cm]
\D \hspace{-0.4cm}& \!\leq \!\!& \!\!(x\!-\!x_i)^{p-1}\!\bigg[\!\frac{(n+2)^{p-1}}{p+1}\!\biggl(\!\!\frac{1}{2^{p-1}n^p}\! +\! \frac{2S^{*}_{p}(n)}{n^{2p+1}}\!\biggr)h_i^{p+1}\|u''\|^p_{\infty} \!+
\!  \frac{1}{3n}\!\bigg(\!\frac{3}{8}\!\bigg)^{\!p}
\!\!h_i^{p+1}\!(\!M_2\!-\!m_2\!)^p\!\bigg]\!\,.
\end{eqnarray*}
It remains now to integrate this inequality on $[x_i,x_{i+1}]$ to obtain that
$$
\D \int_{x_i}^{x_{i+1}}|u(x)-u_{I}(x)|^p dx \leq \frac{(n+2)^{p-1}}{p+1}\biggl(\!\frac{1}{2^{p-1}n^p} + \frac{2S^{*}_{p}(n)}{n^{2p+1}}\!\biggr)\frac{h_i^{2p+1}}{p}\!\|u''\|^p_{\infty}
+ \frac{1}{3n}\bigg(\frac{3}{8}\bigg)^{\!p}
\frac{h_i^{2p+1}}{p}(M_2-m_2)^p\,,
$$
and summing over all values of $i$ from $0$ to $N$ implies that
\begin{equation}\label{N1_u2_V2}
\D \|u-u_{I}\|^p_{0,p} \leq \frac{(n+2)^{p-1}}{p+1}\biggl(\!\frac{1}{2^{p-1}n^p} + \frac{2S^{*}_{p}(n)}{n^{2p+1}}\!\biggr)\frac{h^{2p}}{p}\|u''\|^p_{\infty} +
\frac{1}{3n}\bigg(\frac{3}{8}\bigg)^{\!p}
\frac{h^{2p}}{p}(M_2-m_2)^p.
\end{equation}
Finally, the $W^{1,1}$-norm of the $P_1$-interpolation error is given by adding inequalities (\ref{N1_u1_V2}) and (\ref{N1_u2_V2}) that gives estimate (\ref{Error_Estim_TL}).
\end{prooff}
To compare the error estimates (\ref{Error_Estim_T}) and (\ref{Error_Estim_TL}), we study now the asymptotic behavior of the sum $S^{*}_{p}(n)$ defined in (\ref{Spn}).
For this, we prove the following result.
\begin{lemma}\label{Pascal}
Let $p$ be a non-zero integer and the sum $S_{p}(n)$ defined by: $\D S_{p}(n) = \sum \limits_{k=1}^{n}k^{p}$. \sa
Then, $S_{p}(n)$ have the following asymptotic behavior:
$$
\D S_{p}(n) \sim \frac{n^{p+1}}{p+1}, \, \mbox{ when } \, n \, \mbox{ goes to } + \infty.
$$
\end{lemma}
\begin{prooff}
To prove this result, we proceed by induction on $p$.\sa
First of all, for $p=1$ we have:
$$
\D S_{1}(n) = \sum \limits_{k=1}^{n}k = \frac{n(n+1)}{2}  \underset{n \to +\infty}{\sim} \frac{n^{2}}{2}.
$$
Let us assume  the induction assumption, namely
\begin{equation}\label{Hypo_reccurence}
\D \forall j=1,\dots,p: S_{j}(n) \underset{n \to +\infty}{\sim} \frac{n^{j+1}}{j+1}\,.
\end{equation}
We have to prove that
$$
\D S_{p+1}(n) \underset{n \to +\infty}{\sim} \frac{n^{p+2}}{p+2}.
$$
To this end, recall first  the formula corresponding to a special case of those derived by Blaise Pascal \cite{Pascal}, allowing us to compute $S_{p+1}(n)$ as
\begin{equation}\label{sp_plus_1n}
\D S_{p+1}(n) =\frac{1}{p+2}\bigg[n^{p+2} - 1 - \sum_{j=0}^{p}\binom{p+2}{j}S_j(n)\bigg]\,,
\end{equation}
where we set $\D \binom{p+2}{j} = \frac{(p+2)!}{j!(p+2-j)!}$.\sa
To get the asymptotic behavior of $S_{p+1}(n)$, we notice the next two points:
\begin{enumerate}
\item Due to the binomial expansion, we have:
$$
\D (n+1)^{p+2}-1 = \sum_{k=0}^{p+2}\binom{p+2}{k}n^k - 1 = \sum_{k=1}^{p+2}\binom{p+2}{k}n^k \underset{n \to +\infty}{\sim} n^{p+2}.
$$
\item With the induction assumption (\ref{Hypo_reccurence}), and because $S_j(n)$ is a polynomial of degree less than or equal to $j$, we can write that
$$
\D \sum_{j=0}^{p}\binom{p+2}{j}S_j(n)\, \underset{n \to +\infty}{\sim} \, \sum_{j=0}^{p}\binom{p+2}{j}\frac{n^{j+1}}{j+1} \, \underset{n \to +\infty}{\sim} \, \binom{p+2}{p}\frac{n^{p+1}}{p+1}.
$$
\end{enumerate}
With these two points, it follows from (\ref{sp_plus_1n}) that
$$
\D S_{p+1}(n) \underset{n \to +\infty}{\sim}\, \frac{n^{p+2}}{p+2},
$$
which ends the proof of the lemma.
\end{prooff}
Our goal is now to compare the behavior of the error estimate (\ref{Error_Estim_T}) with the asymptotic one associated with (\ref{Error_Estim_TL}) as $n$ tends to infinity.
Recall that the first one is based on the classical Taylor's formula, whereas the second one is based on Taylor-like formula (\ref{Generalized_Taylor}). \sa
To this end, in the following lemma, we first derive the asymptotic error estimate one can get from (\ref{Error_Estim_TL}), when $n$ tends to infinity.
\begin{lemma}
Let u be in $C^2([0,1])$ and $u_I$ the corresponding $P_1$-interpolation polynomial. Then, the Taylor-like formula (\ref{Generalized_Taylor}) leads to the following asymptotic interpolation error estimate:
\begin{equation}\label{Error_Estim_TL_Asymp}
\forall p\in\N^*: \D \|u-u_{I}\|^p_{1,p} \leq \frac{2}{(p+1)(p+2)}\,h^p\bigg(1 + \frac{h^{p}}{p}\bigg)\|u''\|^p_{\infty}.
\end{equation}
\end{lemma}
\begin{prooff}
The error estimate (\ref{Error_Estim_TL}) being valid for all integer $n\in\N^*$, we are interested in the asymptotic behavior of the error estimate obtained by letting $n$
goes to $+\infty$. \sa
Denoting by $R_n$ the right-hand side of (\ref{Error_Estim_TL}), that is
$$
\D R_n = \frac{(n+2)^{p-1}}{p+1}\biggl(\!\frac{1}{2^{p-1}n^p} + \frac{2S^{*}_{p}(n)}{n^{2p+1}}\!\biggr)\bigg(h^p + \frac{h^{2p}}{p}\bigg)\|u''\|^p_{\infty} + \frac{1}{3n}\bigg(\frac{3}{8}\bigg)^{\!p}\bigg(h^p + \frac{h^{2p}}{p}\bigg)(M_2-m_2)^p\,,\vspace{0.1cm}
$$
it can be decomposed in three parts. \sa
For the first one, we have
$$
\D\lim_{n\rightarrow\infty} \D \frac{(n+2)^{p-1}}{p+1}\frac{1}{2^{p-1}n^p} = \D\lim_{n\rightarrow\infty}
\frac{1}{(p+1)2^{p-1}}\frac{1}{n} = 0\,.
$$
The limit of the second one is obtained with Lemma \ref{Pascal} that leads to
$$
\D\lim_{n\rightarrow\infty} \D \frac{(n+2)^{p-1}}{p+1}\frac{2S^*_p(n)}{n^{2p+1}}
=\D\lim_{n\rightarrow\infty}2\frac{(n+2)^{p-1}}{p+1}\frac{(n-1)^{p+2}}{(p+2)n^{2p+1}}
=\frac{2}{(p+1)(p+2)}\,.
$$
For the third part, we readily get that
$$\D \lim_{n\rightarrow\infty}\frac{1}{3n}\bigg(\frac{3}{8}\bigg)^{\!p}\bigg(h^p + \frac{h^{2p}}{p}\bigg)(M_2-m_2)^p=0.
$$
Putting all together, we obtain that the limit of error estimate (\ref{Error_Estim_TL}) gives (\ref{Error_Estim_TL_Asymp}) when $n\rightarrow + \infty$.
\end{prooff}
\noindent Let us summarize in the following table the constants obtained in the different $W^{1,p}-$ interpolation error estimates, and let us give some examples.
\renewcommand{\arraystretch} {2.25}
\begin{center}
   \begin{tabular}{ | l | c | }
     \hline
     Standard Taylor Theorem & $\D\frac{2^{p-1}}{p+1} +\frac{1}{2}$ \\[0.2cm] \hline
     Mean Value Theorem &  $\D\frac{1}{p+1}$\\[0.2cm] \hline
     Taylor-like Theorem & $\D\frac{2}{(p+2)(p+1)}$ \\[0.2cm]
     \hline
   \end{tabular}
 \end{center}
 \vspace*{0.8cm}

\noindent We notice that the constant of the Taylor-like formula is $2/(p+2)$ smaller than the one obtained with  the mean value theorem. Moreover, standard Taylor formula
leads to a constant strictly greater that 1, since $\D \frac{2^{p-1}}{p+1} \geq \frac{1}{2}$, for $p \geq 1$. \sa
Let us evaluate the improvement obtained in the Taylor-like Theorem by considering a particular case, for instance the Hilbert case $p=2$. In these conditions, the corresponding error estimates are written as
\begin{eqnarray}
\D \mbox{Taylor} & : & \|u-u_{I}\|_{1,2} \leq \bigg(\frac{7}{6}\bigg)^{\!\!\frac{1}{2}}\,h\bigg(1 + \frac{h^{2}}{2}\bigg)^{\!\!\frac{1}{2}}\|u''\|_{\infty}, \label{H1_V1}\\[0.2cm]
\D \mbox{Taylor-like} & : & \|u-u_{I}\|_{1,2} \leq \bigg(\frac{1}{6}\bigg)^{\!\!\frac{1}{2}}\,h\bigg(1 + \frac{h^{2}}{2}\bigg)^{\!\!\frac{1}{2}}\|u''\|_{\infty}, \label{H1_V2}
\end{eqnarray}
so that the constant involved in (\ref{H1_V2}) is $\sqrt{7}$ times smaller compared to the  standard Taylor's formula. \sa
If we consider now the error estimate (\ref{Error_Estim_ACF}) derived by using the mean value theorem, we have
\begin{eqnarray}
\D \mbox{Mean value theorem} & : & \|u-u_{I}\|_{1,2} \leq \bigg(\frac{1}{3}\bigg)^{\!\!\frac{1}{2}} h\bigg(1 + \frac{h^{2}}{2}\bigg)^{\!\!\frac{1}{2}}\|u''\|_{\infty}, \label{H1_V3}
\end{eqnarray}
and the constant obtained by the Taylor-like formula in (\ref{H1_V2}) is still $\sqrt{2}$ times smaller than the one derived in (\ref{H1_V3}).\sa
Let us also illustrate our result with a non-Hilbert case, for example by choosing $p=5$. The corresponding error estimates are expressed by
\begin{eqnarray}
\D \mbox{Taylor} & : & \|u-u_{I}\|_{1,5} \leq \bigg(\frac{19}{6}\bigg)^{\!\!\frac{1}{5}} \,h\bigg(1 + \frac{h^{5}}{5}\bigg)^{\!\!\frac{1}{5}}\|u''\|_{\infty}, \label{L5_1}\\[0.2cm]
\D \mbox{Taylor-like} & : & \|u-u_{I}\|_{1,5} \leq \bigg(\frac{1}{21}\bigg)^{\!\!\frac{1}{5}}\,h\bigg(1 + \frac{h^{5}}{5}\bigg)^{\!\!\frac{1}{5}}\|u''\|_{\infty}, \label{L5_2}
\end{eqnarray}
and the constant is almost 2.5 times smaller in the case of Taylor-like formula.\sa
Here again, considering the mean value theorem, the error estimate is written
\begin{eqnarray}
\D \mbox{Mean value theorem} & : & \|u-u_{I}\|_{1,5} \leq \bigg(\frac{1}{6}\bigg)^{\!\!\frac{1}{5}}\, h\bigg(1 + \frac{h^{5}}{5}\bigg)^{\!\!\frac{1}{5}}\|u''\|_{\infty},\label{New_REF1}
\end{eqnarray}
and the constant obtained in (\ref{L5_2}) is still about 1.3 times smaller than in (\ref{New_REF1}).

\begin{remark}
For the sake of completeness, we can also illustrate the behavior of the Taylor-like error estimate (\ref{Error_Estim_TL}) when $n$ is finite. As above, consider first the Hilbert case $p=2$, and for instance, $n=2$. Substituting these values in expression (\ref{Error_Estim_TL}), we obtain that the second term with $(M_2-m_2)^p$ is negligible before  the first one, so that we can approximately write that
$$
 \|u-u_{I}\|_{1,2} \leq \bigg(\frac{1}{4}\bigg)^{\!\!\frac{1}{2}}\,h\bigg(1 + \frac{h^{2}}{2}\bigg)^{\!\!\frac{1}{2}}\|u''\|_{\infty}\,.
$$
Here, the constant is 0.5, which is approximatively half than the one given by Taylor's formula in (\ref{Error_Estim_T}) which is about 1.08, and slightly smaller than the one derived using the Mean value Theorem in (\ref{Error_Estim_ACF}), which is 0.57.\sa
Let us consider now a non-Hilbert case, also by choosing $p=5$ and $n=2$. The error estimate corresponding to (\ref{Error_Estim_TL}) is written as
$$
 \|u-u_{I}\|_{1,5} \leq \bigg(\frac{1}{8}\bigg)^{\!\!\frac{1}{5}}\,h\bigg(1 + \frac{h^{5}}{5}\bigg)^{\!\!\frac{1}{5}}\|u''\|_{\infty}\,.
$$
The constant $(1/8)^{\frac{1}{5}}\simeq 0.66$ is still two times smaller than those computed by (\ref{Error_Estim_T}) which is about 1.26, and of the same order of magnitude as the one obtained by (\ref{Error_Estim_ACF}), approximatively 0.70.
\end{remark}
\noindent In the next section, we will consider applications of these results to finite element method.
\section{Application to $P_1$ finite element approximation error estimate}\label{Finite_Element}
\noindent The aim of this section is to illustrate, in a simple example,  how we can apply our new results in the context of finite element approximation.\sa
Let $f$ be a given function that belongs to $L^{p}(]0,1[),$ and $u\in W^{2,p}(]0,1[)$ solution to:
\begin{empheq}[left=\mbox{}\hs\empheqlbrace]{alignat=2}
%
\D\hs  -u''(x) + u(x) = f(x), \, x\in\,]0,1[, \hspace{2.5cm}& & & \nonumber \\[0.1cm]
\D  u(0) = u(1)= 0, \hspace{5cm} & & & \nonumber
\end{empheq}
The corresponding variational formulation is given by:
\begin{empheq}[left=\mbox{}\hs\empheqlbrace]{alignat=2}
\mbox{Find } u \in   W^{1,p}_{0}(]0,1[), \mbox{ solution to:} \hspace{5.5cm} & & & \nonumber \\[0.1cm]
\D  \D \int_{0}^{1}\left[ u'(x)v'(x)+u(x)v(x)\right] dx = \int_{0}^{1}f(x)v(x) \,dx , \forall v \in W^{1,q}_{0}(]0,1[), & & & \label{VP_0}
\end{empheq}
where $p$ and $q$ are conjugated, i.e. satisfy (\ref{pq}), and $W^{1,p}_{0}(]0,1[)$ denotes the space of functions $v$ of $W^{1,p}(]0,1[)$ such that $v(0)=v(1)=0$. We
notice that all the integrals in (\ref{VP_0}) are bounded due to Hölder's inequality.\sa
Let us now introduce the finite-dimensional subspace $V_h$ of $W^{1,q}_{0}(]0,1[)$, consisting of functions $v_h$ defined on $[0,1]$ which are piecewise linear on each interval $[x_i,x_{i+1}], (i=0,N)$, and satisfying the boundary conditions: $v_h(0)=v_h(1)=0$. We also consider $u_h$ the approximation of the solution $u$.\sa
To apply error estimates (\ref{Error_Estim_T}) and (\ref{Error_Estim_TL}), we also assume that solution $u$ belongs to $C^2([0,1])$, which is consistent with the regularity of $W^{1,p}(]0,1[)$, as well as with the regularity of the finite-dimensional subspace $V_h$ defined above. \sa
As a first example, let us consider the Hilbert case, i.e. with $p=q=2$, and $W^{1,p}(]0,1[)=H^1(]0,1[$. In this case, we can apply the classical Céa's Lemma
\cite{Ern_Guermond}, \cite{BrSc08}, which states that for all $v_{h} \in V_{h}$,
\begin{equation}\label{CEA_Hilbert}
\left\|u - u_{h} \right\|_{H^1} \leq C \|u-v_h\|_{H^1},
\end{equation}
where $C$ is a positive constant which depends on the bilinear form introduced in (\ref{VP_0}). \sa
Then, in (\ref{CEA_Hilbert}) we choose the particular function $v_h$ defined by $v_h=u_I$, where $u_I$ denotes the interpolation polynomial of the solution $u$ satisfying
the boundary conditions $u_I(0)=u_I(1)=0$. \sa
Therefore, (\ref{CEA_Hilbert}) leads to
\begin{equation}\label{CEA_u_I}
\left\|u - u_{h} \right\|_{H^1} \leq C \|u-u_I\|_{H^1}\,.
\end{equation}
Now, let us assume a first mesh with a mesh size equal to $h_1$, where the classical Taylor method is used, resulting in the error estimate (\ref{H1_V1}). Also, assume a second mesh of mesh size $h_2$ which is concerned by the  Taylor-like estimate (\ref{H1_V2}). Hence, by (\ref{CEA_Hilbert}), we readily see that
\begin{eqnarray*}
\D \mbox{Taylor} & : & \left\|u - u_{h} \right\|_{H^1} \leq C\,\bigg(\frac{7}{6}\bigg)^{\!\!\frac{1}{2}}\,h_1\bigg(1 + \frac{h_1^{2}}{2}\bigg)^{\!\!\frac{1}{2}}\|u''\|_{\infty} \simeq C\,\bigg(\frac{7}{6}\bigg)^{\!\!\frac{1}{2}}\, h_1\,\,\|u''\|_{\infty}, \\[0.2cm]
\D \mbox{Taylor-like} & : & \left\|u - u_{h} \right\|_{H^1} \leq C\,\bigg(\frac{1}{6}\bigg)^{\!\!\frac{1}{2}}\,h_2 \bigg(1 + \frac{h_2^{2}}{2}\bigg)^{\!\!\frac{1}{2}}
\,\|u''\|_{\infty} \simeq C\,\bigg(\frac{7}{6}\bigg)^{\!\!\frac{1}{2}}\, h_2\,\|u''\|_{\infty}.
\end{eqnarray*}
Consequently, if we want to ensure that the approximation error $\left\|u - u_{h} \right\|_{H^1}$ is smaller than a specified threshold, the above estimates leads to: $h_2 = \sqrt{7} h_1$.\sa
In other words, $h_2$ may be around 2.65 times greater than $h_1$.  A practical consequence is the possibility of using a coarser mesh for a given accuracy. This reduction in
terms of the total numbers of meshes would be even more significant when considering the extension of this case to three-dimensional applications. Indeed, assuming that the
three axis are discretized similarly, the number of nodes could be about $2.65^{3}\simeq 18$ times fewer,   significantly reducing the cost of finite element implementation.\sa
For the second example, we consider a non-Hilbert case. We take $p=5$ and $q=5/4$ and apply the generalized C\'ea's Lemma \cite{Ern_Guermond} valid in Banach spaces, which
is expressed as follows:
$$
\D\|u-u_h\|_{W} \leq \left(1+\frac{\|a\|_{W,V}}{\alpha_h}\right)\inf_{w_h\in W_h}\|u-w_h\|_{W} \,.
$$
In this lemma, two  different (conjugated) spaces are involved. In our case, $W=W^{1,5}_{0}(]0,1[)$ and $V=W^{1,5/4}_{0}(]0,1[)$, which are a Banach space and a reflexive Banach space, respectively, as required.
\sa
Under these conditions, we readily obtain from this lemma the following inequality:
\begin{equation}\label{Cea_W1p}
\D\|u-u_h\|_{W^{1,5}}
\leq \left(1+\frac{\|a\|_{W^{1,5},W^{1,5/4}}}{\alpha_h}\right)\|u-u_I\|_{W^{1,5}},
\end{equation}
where $u_I$ here also denotes the interpolation polynomial of the solution $u$.\sa
Now, consider two given meshes characterized by their mesh sizes $h_1$ and $h_2$ used under the same conditions as above. \sa
Then, if we assume than $W_h=V_h$, from estimate (\ref{Cea_W1p}) and (\ref{L5_1})-(\ref{L5_2}), we obtain the  following two approximation error estimates
$$
\hspace*{-0.45cm}\mbox{Taylor}  :  \left\|u - u_{h} \right\|_{1,5} \leq  C'\,\left\|u - u_I \right\|_{1,5} \leq C'\,\bigg(\frac{19}{6}\bigg)^{\!\!\frac{1}{5}} \, h_1 \,\bigg(1 + \frac{h_1^5}{5}\bigg)^{\!\!\frac{1}{5}} \, \|u''\|_{\infty} \simeq C'\,\bigg(\frac{19}{6}\bigg)^{\!\!\frac{1}{5}}\,h_1\,\|u''\|_{\infty},
$$
$$
\hspace*{0.1cm}\D \mbox{Taylor-like} :  \left\|u - u_{h} \right\|_{1,5} \leq  C'\,\left\|u - u_I \right\|_{1,5} \leq C'\,\bigg(\frac{1}{21}\bigg)^{\!\!\frac{1}{5}}\,h_2\, \bigg(1 + \frac{h_2^5}{5}\bigg)^{\!\!\frac{1}{5}}\, \|u''\|_{\infty} \simeq C'\,\bigg(\frac{1}{21}\bigg)^{\!\!\frac{1}{5}}\,h_2\,\|u''\|_{\infty},
$$
where we set: $\D C' = 1+\frac{\|a\|_{W^{1,5},W^{1,5/4}}}{\alpha_h}$.\sa
Consequently, if we want to ensure that the upper bound of the approximation error $\left\|u - u_{h} \right\|_{1,5}$ is smaller than a specified threshold, we find from the above estimates that $h_2\simeq 2.31 \,h_1$. Consequently, the reduction in the number of nodes is approximately $2.31^3\simeq 12.4$ times less for 3D applications.

\section{Conclusion}\label{Conclusion}
\noindent In this paper we derived several $W^{1,p}$ $P_1$-interpolation error estimates. We obtained them using first the standard Taylor's formula, then based on the mean value theorem, and finally, by a Taylor-like formula. \sa
These estimates were derived by applying Fubini's theorem and Hölder's inequality while maintaining a consistent approach throughout. This guarantees a unified methodology enabling us to compare the various constants obtained in theses estimates.
\sa
These results allow us to emphasize that the use of the Taylor-like formula leads to a significant gain, as the corresponding constant in error estimate is smaller than those associated with the standard Taylor's formula or the mean value theorem. In particular, we highlighted that the constant that appears with Taylor's formula is strictly greater than one, whereas mean value theorem leads to a constant  smaller than one. The Taylor-like formula, for its part, gives a constant that is $2/(p+2)$ smaller than the one derived from the mean value theorem.
\sa
We also illustrate our results in the context of the finite element method. As an example, we introduce a
second-order differential equation with a right-hand side in $L^{p}(]0,1[)$. We then consider the Hilbert case where $p=2$ and its corresponding error estimate in $H^1(]0,1[)$, as well as the Banach case when $p=5$ along with its error estimate in $W^{1,5}(]0,1[)$. In this way, we showed that the number of nodes needed for a given mesh is about 18 times less when $p=2$ and 12 times less when $p=5$, assuming a three-dimensional mesh.\sa
The outlook of this research essentially involves extending these estimates to the case where the dimension of the domain of interpolation, and then of integration, is strictly greater than one. This will require a Taylor-like formula we already derived in \cite{ChAs2023}. On the other hand, we will also explore an extension to a second-order Taylor-like formula, as we proved in \cite{JCAM2023}. Both extensions will be examined to evaluate their impact in the error estimates in the context of numerical analysis applications.
\sa
\noindent \textbf{\underline{Homages}:} The authors want to warmly dedicate this research to pay homage to the memory of Professors Andr\'e Avez and G\'erard Tronel who largely promote the passion of research and teaching in mathematics of their students.

\end{document}